# Mathematical Knowledge and the Role of an Observer

## Ontological and epistemological aspects

Mark Burgin

*University of California, Los Angeles*
Los Angeles, CA 90095, USA

**Abstract.** As David Berlinski writes (1997), "*The existence and nature of mathematics is a more compelling and far deeper problem than any of the problems raised by mathematics itself.*" Here we analyze the essence of mathematics making the main emphasis on mathematics as an advanced system of knowledge. This knowledge consists of structures and represents structures, existence of which depends on observers in a nonstandard way. Structural nature of mathematics explains its reasonable effectiveness.

1. Introduction

This paper is inspired by the article of Yu.I. Manin "Mathematical Knowledge: Internal, social and cultural aspects" (2007). In his article, Manin discusses the essence of mathematical knowledge, the role of a mathematics in human culture and the place of mathematics in society. At the beginning of his interesting article, he proclaims that "Pure mathematics is an immense organism built entirely and exclusively of ideas that emerge in the minds of mathematicians and live within these minds."

After this, Manin suggests three possibilities of other images of mathematics:

◊ *Mathematics is the contents of mathematical manuscripts, books, papers and lectures*, with the increasingly growing net of theorems, definitions, proofs, constructions, conjectures, i.e., the results of the activity of mathematicians.

◊ *Mathematics is the activity of mathematicians*.

◊ *Mathematics is a part of the Platonic world of ideas*.

Putting forward here three other but related perspectives on mathematics, we base our approach on the large-scale structure of the world has the form of the *Existential Triad of the World* (Burgin, 2012), which is presented in Figure 1.

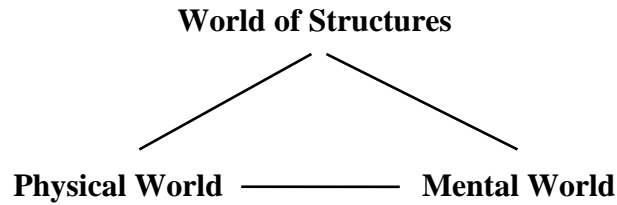

**Figure 1.** The Existential Triad of the World

In the Existential Triad, the Physical (material) World is conceived as the physical reality studied by natural sciences, the Mental World encompasses different levels of mental reality, the lowest levels of which are studied by psychology, and the World of Structures consists of various forms and types of structures, an essential part of which is studied by mathematics.

The existential stratification of the World continues the tradition of Plato with his World of Ideas (Plato, 1961) and the tradition of Charles Sanders Peirce with his extensive triadic classifications (Peirce, 1931-1935).

An interesting peculiarity of our days is that if before people reduced the whole world to its physical component, now there are suggestions to reduce the whole world to its structural component (cf., for example, (Tegmark, 2008)).

The Existential Triad brings into play three interrelated views on mathematics:

- *Mathematics as mathematical knowledge* (the structural perspective).
- *Mathematics as a social system* (the material perspective).
- *Mathematics as a part of people's mentality* (the mental perspective).

These three interpretations of mathematics do not exclude or contradict one another. To the contrary, they complement one another reflecting projections of mathematics on the components of the Existential Triad. As a result, only this unity of interpretations allows achieving a comprehensive understanding of mathematics.

Note that writing about people's mentality, we have in mind not only individual mentality but also group mentality and social mentality, which constitute the second and the third levels of the Mental World.

Here we extend and specify these visions of mathematics beginning with an analysis of mathematical knowledge treating it as knowledge about mathematical objects.

## 2. Mathematics as the science of structures

The ontological analysis of mathematical knowledge demands answering two fundamental questions – what objects mathematics studies and where these objects are situated. We analyze the first question in this section and deal with the second one in the fourth section.

Discussing objects of mathematical knowledge, Manin follows Paul Davis and Reuben Hersh defining mathematical objects as "ideas which can be handled as if they were real things" although he notes that Davis and Hersh call them "mental objects with reproducible properties."

However, this definition does not specify mathematical knowledge because, for example, philosophy also studies ideas. Even more, it is possible to argue that natural sciences, such as physics, also study ideas, namely, ideas about nature, while social sciences study ideas about society. Thus, it is reasonable to look for another interpretation and/or understanding of mathematical objects.

To figure out the scope and domain of mathematics, it is useful to look at the history of mathematics. At the beginning, mathematics was the discipline that studied numbers and geometrical shapes/forms (Burton, 1997). However, the development of mathematics changed its essence and mathematics has become the discipline that studies exclusively mathematical structures because taking any mathematical object, we can easily see that it is a structure of some kind. Formal structures occupy the central place in the unified picture of mathematics. In this picture, numbers and geometrical shapes/forms, which are still actively studied by mathematicians, are and are also treated as definite kinds of structures. For instance, according to Shapiro (1997), "the essence of a natural number is the relations it has with other natural numbers."

Following Lautman (1938), the most active proponent of this approach Bourbaki defines pure mathematics as nothing other than the study of pure structure (Bourbaki, 1948; 1957; 1960). In particular, Bourbaki (although some ascribe these words exclusively to Dieudonné) describes the unifying role of the structures in the following way:

"Each structure carries with it its own language, freighted with special intuitive references derived from the theories from which the axiomatic analysis ... has derived the structure. And, for the research worker who suddenly discovers this structure in the phenomena which he is studying, it is like a sudden modulation which orients at once the stroke in an unexpected direction in the intuitive course of his thought and which illumines with a new light the mathematical landscape in which he is moving about.... Mathematics has less than ever been reduced to a purely mechanical game of isolated formulas; more than ever does intuition dominate in the genesis of discoveries. But henceforth, it possesses the powerful tools furnished by the theory of the great types of structures; in a single view, it sweeps over immense domains, now unified by the axiomatic method, but which were formerly in a completely chaotic state." (Bourbaki, 1950)

Some philosophers and many mathematicians also assume that all mathematical objects are structures of some type. For instance, Resnick (1997) explains:

"In mathematics … we do not have objects with an 'internal' composition arranged in structures, we have only structures. The objects of mathematics … are structureless points or positions in structures. As positions in structures, they have no identity or features outside of a structure."

It is interesting that treating structures in mathematics, Bourbaki define and discuss inner structures, while Resnick talks about outer structures in the sense of the general theory of structures (Burgin, 2012).

As North (2009) writes, the idea of mathematical structure is relatively straightforward. There are number structures, geometric structures, topological structures, algebraic structures, and so forth. Mathematical structure tells us how simpler mathematical objects are organized to form more complex mathematical objects. Different sources provide various definitions of mathematical structures given by different authors while the most comprehensive definition of a structure in general is developed in the general theory of structures (Burgin, 2012).

Here we do not give this definition but describe the basic types of structures discovered in the general theory of structures.

First, there are abstract structures and embodied structures.

An *embodied structure Q* is assigned to another system $R$ and treated as a structure of $R$. In this case, it is possible to say that $Q$ is metaphorically embodied in $R$.

An *abstract structure Q* is not assigned to any other system and can be treated only as a structure of itself.

Note that an embodied structure can be a structure of a physical, mental or structural system. For instance, a three-dimensional sphere is a structure of the Earth (a material system), of the mental image of the Earth (a mental system), and of a section of a four-dimensional sphere (a structural system).

Second, according to the general theory of structures (Burgin, 2012), there are five types of embodied structures:

1. An *internal structure TQ* of a system *R* contains only inner structural parts, components and elements, i.e., parts, components and elements of *R*, relations between these parts, components and elements, relations between these parts, components, elements and relations from *TQ* and relations between relations from *TQ*.
2. An *inner structure IQ* of a system *R* is a substructure of an internal structure *TQ* of *R*, where *IQ* is obtained by exclusion of (1) the whole system *R* as a part, component or element of itself and (2) all relations that include *R*.
3. An *external structure EQ* of a system *R* is an extension of the internal structure, in which other systems, their parts, components and elements are included, as well as relations between all these included parts, components and elements, relations between these parts, components, elements and relations from *EQ* and relations between relations from *EQ*.
4. An *intermediate structure MQ* of a system *R* is a substructure of an external structure *EQ* of *R*, where *MQ* is obtained by exclusion of (1) the whole system *R* and other systems from *EQ*, as well as (2) all relations that include these systems.
5. An *outer structure OQ* of a system *R* is an inner structure of a system *U* in which *R* is only one of the inner elements of the inner structure *IQ* of the system *U*.

All structures as abstract objects belong to the World of Structures (cf. Figure 1). At the same time, all structures embodied in physical systems belong to the Physical World, while all structures embodied in mental systems belong to the Mental World.

Note that all fields of mathematics study specific structures: geometry studies geometrical shapes, algebra studies algebraic systems such as groups, rings and fields, and functional analysis studies functions, functionals and operators. In contrast to this, the general theory of structures studies structures *per se* and contains three parts: philosophical, methodological and mathematical.

To better understand structure embodiment, we need to know that a structure embodied in a physical system, which is thus a structure of this system, can have not only physical elements but also ideal elements. For instance, taking a circle drawn on a board or even in sand, we know that elements of this physical circle are not points but little particles, e.g., grains of sand in the second case, but one of the structures of this physical circle is a mathematical circle, which consists of ideal points. Besides, the physical circle is discrete while its structure is a mathematical circle, which is continuous.

It is necessary to remark that not all mathematicians assume that all mathematical objects, and in particular, the set of natural numbers, are structures. For instance, Rodin (2011) writes about categories without structures. However, this is point of view is a consequence of a very limited understanding of the concept *structure*. An adequate understanding tells us that a category is a structure and thus, the title of Rodin's paper can be read as "Structures without structures".

Besides, the diversity of interpretations and applications of mathematical objects brought researchers to the idea that mathematics is a science of patterns (cf., for example, (Steen, 1988; Resnick, 1997)). As Steen (1988) writes, the rapid growth of computing and applications has helped cross-fertilize the mathematical sciences, yielding an unprecedented abundance of new methods, theories, and models. Examples from statistical science, core mathematics, and applied mathematics illustrate these changes, which have both broadened and enriched the relation between mathematics and science. No longer just the study of number and space, mathematical science has become the science of patterns, with theory built on relations among patterns and on applications derived from the fit between pattern and observation. However, there are individual and pattern mathematical structures and mathematics studies both kinds.

To better delineate the content of mathematics, it is necessary to understand that there are many structures that are not mathematical. For instance, there is a variety of chemical structures. Moreover, structural realism assumes that science studies structures of systems in nature and can give knowledge only about these structures but not about natural systems as they are. Thus, mathematics studies an essential part of the World of Structures but not everything that exists in this World.

Mathematicians study structures as biologists study living organisms or astronomers study celestial bodies. As other scientist, mathematicians make observations of the studied reality, i.e., of known structures; perform measurements in structural reality calculating and logically inferring properties of studied structures; devise narrowly targeted experiments for finding and proving statements about properties of studied structures; and finally produce a mathematical law (theorem) or a theory, which becomes a current milestone of mathematics (Burgin, 1998).

Naturally, we come to the problem of existence of mathematical and other structures. To rigorously treat this problem, at first, we explore the meaning of existence of physical (material) objects (things) in the Physical World.

### 3. Existence of material things and the role of an observer

What does it mean that a material thing, such as a mountain, table or chair, exists?

There is a naïve point of view that if people can comprehend something by their senses, then it exists. However, attentive people were able to see that there are situations when senses deceive creating only an illusion of real existence. As a reflection of this experience, ancient Indians created the concept of Maya, which can be interpreted as illusion, which substitutes reality, as well as the power to create such an illusion.

In Europe, Plato introduced the idea that the material world is only a reflection of the real World of Ideas. This approach also undermined absolute reality of the material world.

However, the majority of people believed in real existence of material things they could comprehend by their senses, e.g., see, hear, etc. This brought philosophers to the problem of an observer. The outstanding philosopher George Berkeley was, may be, the first to formulate it in 1710 writing:

"But, say you, surely there is nothing easier than for me to imagine trees, for instance, in a park [...] and nobody by to perceive them. [...] The objects of sense exist only when they are perceived; the trees therefore are in the garden [...] no longer than while there is somebody by to perceive them." (Berkeley, 1710)

The problem of the role of an observer in physical reality became very urgent for physicists after the quantum world had been discovered and studied. The reason was that existence in the quantum world is very different from existence of things in our everyday life because to find that something exists in the quantum, it is necessary to organize sophisticated experiments, which are explicitly or implicitly based on advanced theories. This vague situation with quantum brought back Berkeley's consideration. For instance, Albert Einstein once asked Niels Bohr if the Moon exists when no one is looking at it. Another questions asked by physicists were whether there is sound in a forest when a tree falls if there is no one there to observe it or whether a tree has fallen if there is no observer.

Beedham explains this situation in the following way.

"It is indeed an opinion strangely prevailing amongst men, that houses, mountains, rivers, and in a word all sensible objects have an existence natural or real, distinct from their being perceived by the understanding. […]

The notion that without an observer the world does not exist is counterintuitive, and appears to contradict our immediate apprehension of the world. But many things in life – scientific facts, as they say – are counter-intuitive.

We used to think the earth was flat: it took physics, intrepid travellers, and photographs from outer space to show us that the earth is round. And yet it still seems bizarre and contrary to our immediate experience that we are stuck to a ball of rock spinning through space, and don't fall off. It sounds equally perverse to claim that the world is a product of the perceptual apparatus of the observer, and that the world does not exist but for the observer, yet it is true." (Beedham, 2005).

The approach advocated in this paper gives answers to all such questions. Contemplating the Einstein's question, it is possible to suggest that the Moon have existed as the Moon only after it was identified as the Moon. Before (for example, the first second after emergence of our universe) there was not such a physical object. Later when the Moon was formed, it acquired a potential existence and only when people identified it as the Moon, its existence became actual.

Let us ask one more question. Does the Mount Everest exist? The answer will be yes for people living in England or USA but it will be no for people living in Tibet, where the same mount is called Chomolungma ("goddess mother of the world") and no for people living in Nepal, where the same mount is called Sagarmatha ("goddess of the sky").

We may also ask if the mount Everest exists as a mount (mountain). Analyzing the situation, we see that without people, it exists no more than as a part of nature. Only when people structure the nature in a definite way, it starts existing as a mount (mountain). For instance, we even do not know if animals or birds perceive it as a specific object or merely as an unspecified part of their surroundings.

The situation becomes more transparent if we ask the question whether the planet Pluto exists. The answer will look a little bit strange for those who believe that material things definitely exist. Namely, the planet Pluto existed from 1930, when it was discovered and (!) named, and until 2006, when it was announced that the celestial body called Pluto was not a planet but a *dwarf planet*. Note that in astronomy, a dwarf planet is not a planet. Likewise, in named set theory, a named set is not always a set (Burgin, 2011).

Thus, before it was discovered and named, the Pluto had potential existence both as a planet and as a dwarf planet. After it was discovered and named, the Pluto had at first actual existence as a planet and then as a dwarf planet.

As a result, we come to the following conclusion. Without an observer, physical things have only potential existence. When an observer comprehends a part of the nature as some object, it acquires actual existence. Note that the task of an observer is not only to indicate existence of something but also to form an image of the comprehended entity.

Here is one more interesting question. Does the mountain in a picture exist? The mountain in a picture is a material object because it is a part of the picture, which is a material object, and according to the general understanding, any material object exists. However, the mountain in a picture is not a real mountain. It exists only as an image of a real mountain but not as a real mountain. This shows that when we speak about existence, we always have in mind existence of some object.

Then again, speaking or thinking about physical (material) existence of an object *A*, people assume that some material thing, which they call *A*, has definite properties that make it *A*. People ascribe these properties to this object and thus, the object does not exist as *A* without people as its observers.

## 4. Existence of structures in general and mathematical structures in particular: Three modes and three worlds

The role of the observer becomes even more vital when we are dealing with the World of Structures (Plato Ideas) in general and mathematical structures, in particular.

In Section 2, analyzing the structure of the world given by the Existential Triad, we came to the conclusion that the most consistent view is that mathematical knowledge is a domain in the World of Structures, elements of which are essentially connected to two other worlds – the physical world and mental world. There are two types of these connections:

- *Representation* of knowledge by physical and mental structures
- *Interpretation* of knowledge in physical and mental structures

In the physical world, mathematical knowledge is represented by mathematical books, journals, papers, drawn figures, pictures, videos, etc.

In the mental world, mathematical knowledge is represented by mathematical ideas, images, concepts, etc.

In the physical world, mathematical knowledge is interpreted by physical systems and processes. For instance, Newton equations are interpreted as mathematical descriptions of dynamics of systems of rigid bodies, e.g., as dynamics of the Solar System.

In the mental world, mathematical knowledge is interpreted by mental systems and processes.

It is interesting that while both representations – mental and physical - of mathematical knowledge are dynamic being in the permanent process of development and emergence, mathematical knowledge per se in the world of structures is invariable in its potentiality.

To understand that structures exist in reality, it is useful to discern different types of reality and three kinds of existence:

*Actual existence* of an object means that this object belongs to the corresponding reality in agreement with what people understand when they speak about such an object.

*Potential existence* of an object means that this object can potentially but does not now belong to the corresponding reality in agreement with what people understand when they speak about such an object.

*Emerging existence* of an object means transformation of this object's state from potential to actual existence.

According to contemporary science, an object is real if it is possible to test its existence by observation and experiment. This is how people substantiate reality of trees, mountains and tables. This is how physicists assert reality of atoms, electromagnetic waves and radiation. This is how astronomers verify reality of planets and stars.

Thinking about mathematical structures, we can easily find examples of similar observations. For instance, looking at a ball, people can without difficulty observe the geometrical structure called a sphere. Looking at a coin, people can effortlessly observe the geometrical structure called a circle. Looking at an envelope, people can certainly observe the geometrical structure called a rectangle and so on and so forth.

In some cases, pure observation is insufficient and we need experiments to observe mathematical structures. For instance, taking some object and performing such an experiment as counting them, people can observe the structure called natural numbers. Thus, the question whether numbers existed before this experiment is the same whether the Moon existed before somebody observed it. It means that existence of (mathematical) structures depends on the observer in the same way as existence of material things.

It is possible to observe not only mathematical structures that came from the everyday experience of people but also even the most fundamental ones. To show this, let us consider the basic structure in mathematics called a *fundamental triad* (or a *named set*) (Burgin, 1990; 2004; 2011), the theory of which form the unified foundation for mathematics, and describe observations that confirm existence of this structure. Here we regard only basic fundamental triads (basic named sets) as the simplest structures, which are still observable.

In the symbolic representation, a fundamental triad (named set) has the form $\mathbf{X} = (X, f, I)$ in which $X$ is called the *support* of $\mathbf{X}$, $I$ is the *component of names* (*reflector*) of $\mathbf{X}$, and $f$ is the *naming correspondence* (*reflection*) of $\mathbf{X}$.

In the graphic representation, the following diagrams describe a fundamental triad:

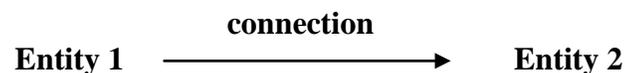

or

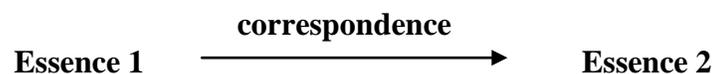

Thus, people can see structures but as Plato wrote about Ideas (eidos), they can comprehend structures only using their intelligence and knowledge.

Now let us describe how we can see a fundamental triad. Imagine, it is raining and you can see a cloud, the ground and streams of water, which go from the cloud to the ground. In essence, it is a fundamental triad $(X, f, I)$, in which the support $X$ is the cloud, the reflector $I$ is the ground and the streams of water build the reflection $f$.

Some will object that they do not see a fundamental triad. They will argue that they can see only a cloud, rain and ground. To understand why it happens, let us look into the history of human civilization.

For a long time, the majority of people have assumed that only those things (objects) were real that were comprehensible by the senses people had. However, even a material thing such as a table, tree or mountain, is comprehended not only by senses alone – understanding is necessary. For instance, if a person who lived in ancient Greece would see a plane, this person would not understand what he saw without explanation, and for many even a good explanation would not be enough.

To clearly see the reality of structures, it is useful to know that in some areas, structures have played an indispensable role for a long time. For instance, in chemistry, researchers have recognized the significance of chemical structures almost from the very beginning of chemistry as a science. May be the most evident example of the importance of structures is the striking difference between diamond and graphite. Both materials (substances) are built entirely from the chemical element carbon and the difference in their properties is completely caused by their distinct crystal structures.

One more transparent example of a structure is an organization because organization is not so much people that work in the organization or material things used by these people but a definite structure. Without the organizational structure, it will be only a bunch of people and things. The same people can form different organizations depending on different organizational structures.

Embedded structures exist in diverse physical and mental things. At the same time, to figure out where the world of abstract structures exists, it is necessary to understand that there are different types of existence and physical or material existence is only one of them. This understanding demands a mental effort similar to the effort in seeing that many stars, which people can see only as points, are actually as big as the Sun and some are even bigger.

Besides, when we say that the Earth existed long before the first living creatures appeared, it is only a reasonable extrapolation of our knowledge. However, as a reasonable extrapolation, we can also assume that the set **N** of all natural numbers had existed before people appeared.

Some philosophers invented the following objection to mathematical Platonism (cf., for example, (Leifer, 2016)). If people have intuitive access to the Platonic realm, our physical brains must interact with it in some way. Our best scientific theories contain no such interaction.

However, the latter statement is not true because scientific theories provide mathematical tools to comprehend (observe) the World of Structures as a scientific incarnation of the World of Ideas. Moreover, mathematically equipped brain can have images of diverse mathematical structures being in such a way connected to the World of Structures.

In more details, existence of the World of Structures and its relation to Plato's World of Ideas is considered in (Burgin, 2012; 2017)

## 5. Conclusions

Mathematics studies mathematical structures as astronomy studies celestial bodies and these structures are real as those celestial bodies are real. Only mathematical structures exist not in the Physical World but in the World of Structures and as Plato suggested, it is necessary to make a mental effort to see this World of Structures.

If we believe that material things are real, similar arguments imply reality of mathematical structures. Some of them are abstract, while others are embodied in material things.

Note that existence of embodied mathematical structures solves the problem formulated by Albert Einstein in 1921 and later reiterated by different authors about "unreasonable" effectiveness of mathematics Eugene Wigner discussed this problem with respect to physics in 1959/1960, Richard Hamming did this with respect to engineering in 1980 and Arthur Lesk did this with respect to molecular biology in 2000. Here we come to the following answer. It is natural that when mathematics correctly describes and utilizes key embodied structures, it reflects the essence of things.

However, existence of mathematical structures embodied in physical bodies does not mean that mathematics studies these bodies and is, as Leifer (2016) suggests, a natural science. Mathematics is a science but it studies the realm of the World of Structures because it studies not materials things but structures in general and structures of materials things in particular (Burgin, 1998).

It is important to understand that existence of structures is akin existence of material bodies. Structures in general and mathematical structures in particular exist not only in minds of people. However, their comprehension demands mental work for constructing mental representations and understanding of structures. This peculiarity was exactly described by Plato if Ideas were interpreted as structures.

Although some think that there is no time in mathematics, at the social and individual levels, mathematical knowledge exists as a dynamical system. Usually it grows in society but there were periods when some knowledge was lost and the whole system of mathematical knowledge declined. For instance, an essential part of the mathematical knowledge of ancient Greece was lost in the medieval Europe.

Similar situations frequently exist on the individual level when people learn something but later forget it. However, social and psychological issues of mathematical knowledge are studied in other works of the author.